\documentclass[11pt,a4paper]{article}

\usepackage{amsmath,amsthm,amsfonts,amssymb,color}
  \definecolor{darkgreen}{rgb}{0,0.4,0}

\usepackage{tikz-cd}

\usepackage{url}

\usepackage[normalem]{ulem}

\parskip=\medskipamount
\def\onehalf{{\textstyle\frac12}}

\def\lie#1{{\mathcal L}_{#1}}
\def\fpd#1#2{\frac{\partial #1}{\partial #2}}
\def\R{{\mathbb R}}
\def\SODE{{\textsc{sode}}}

\def\hook{{\mathchoice{\vrule height 0pt depth 0.4pt width 3pt
\vrule height 5pt depth 0.4pt \kern 3pt} {\vrule height 0pt depth
0.4pt width 3pt \vrule height 5pt depth 0.4pt \kern 3pt} {\vrule
height 0pt depth 0.2pt width 1.5pt \vrule height 3pt depth 0.2pt
width 0.2pt \kern 1pt} {\vrule height 0pt depth 0.2pt width 1.5pt
\vrule height 3pt depth 0.2pt width 0.2pt \kern 1pt} }}

\def\d{\mbox{d}}

\theoremstyle{plain}
\newtheorem{thm}{Theorem}[section]

\newtheorem{propn}[thm]{Proposition}
\newtheorem{cor}[thm]{Corollary}

\theoremstyle{definition}

\newtheorem{xmpl}{Example}

\newtheorem*{xmpl*}{Example}
\newtheorem*{appln*}{Application}

\newcommand{\pd}[2]{\frac{\partial#1}{\partial#2}}
\newcommand{\vf}[1]{\frac{\partial}{\partial #1}}
\newcommand{\bnabla}{\bar{\nabla}}

\newcommand{\Pnabla}{{\nabla}^P}
\newcommand{\Pbnabla}{{\bar\nabla}^P}

\newcommand{\PH}{P_H}

\newcommand{\Sp}{\text{Sp}} 
\newcommand{\Img}{\text{im}}
\newcommand{\Ker}{\text{ker}}
\newcommand{\X}{\mathfrak{X}}
\newcommand{\F}{\mathfrak{F}}
\newcommand{\M}{\mathcal{M}}

\newcommand{\vfe}[2]{\fpd{}{#1}\biggr|_{#2}}

\begin{document}

\title{Covariant derivatives for Ehresmann connections}

\author{G.E.\,Prince and D.J.\,Saunders}

\maketitle

\begin{quote}

{\bf Abstract.} {\small  We deal with the construction of covariant derivatives for some quite general Ehresmann connections on fibre bundles. We show how the introduction of a vertical endomorphism allows construction of covariant derivatives separately on both the vertical and horizontal distributions of the connection which can then be glued together on the total space. We give applications across an important class of tangent bundle cases, frame bundles and the Hopf bundle.\\


}

\end{quote}

\section{Introduction}

We address the question of the existence of covariant derivatives for Ehresmann, or `nonlinear', connections~\cite{E50}. We will show how quite simple additional requirements can yield the familiar covariant derivative of the linear case.

We shall take as our definition of an Ehresmann connection the minimalist version, using only a horizontal distribution on the tangent bundle of a fibre bundle complementary to the vertical sub-bundle:\ although other authors sometimes impose additional conditions involving completeness~(\cite{GHV72}, see also~\cite{KMS93}), we shall not need them.

Using a very recent construction~\cite{PPSD21} of a covariant derivative on a manifold from covariant derivatives on vector submodules or distributions, we show how to construct covariant derivatives on distributions of common rank using certain endomorphisms. In the context of Ehresmann connections on a fibre bundle $\pi:E\to B,$ we then show how to construct covariant derivatives on the corresponding vertical and horizontal distributions (of possibly different ranks) in quite general cases, not necessarily principal or vector bundles, glueing them together to give a covariant derivative on $E$.

This framework finds application in general cases such as  the Hopf bundle; it gives a new understanding of affine, nonlinear and spray connections in the tangent bundle case; and it provides natural covariant derivatives on frame bundles using a novel decomposition of the corresponding vertical sub-bundle.

We remark that a method of constructing a covariant derivative from an Ehresmann connection on a general fibred manifold was given in~\cite{CI87}, but in that approach the operator was defined only on projectable vector fields, mapping them to vertical vector fields.

\section{The construction of covariant derivatives I}
In proposition \ref{propn1} and theorem \ref{thm1} we reprise a new construction of a covariant derivative on a manifold using covariant derivatives on submodules of vector fields, taken from the recent paper \cite{PPSD21} where the proofs can be found in section 2. We begin by demonstrating, as a first step, an extension of these submodule derivatives. In what follows we use the terms submodule and distribution interchangeably even though they are distinct constructs, by considering for any distribution the submodule of vector fields taking their values in that distribution.

Let $\M$ be a manifold and $P$ a projector onto a given distribution.
Suppose that $\Img(P)$ has a covariant derivative $\Pnabla$, that is $\Pnabla_XY$ has the usual properties for $X,Y\in \Img(P)$ over the smooth functions, $\mathfrak{F}(\M),$ on $\M$, namely $\R$-linearity in both arguments and
\begin{equation}\label{cov deriv props}
\Pnabla_{fX}Y=f\Pnabla_XY \ \text{and}\ \Pnabla_X(fY)=X(f)Y+f\Pnabla_XY
\end{equation}
(we do not require that $\Pnabla_XY \in \Img(P)$).
Then $\Pnabla$ can be partially extended to the module of vector fields $\X(\M)$ as follows:

\begin{propn}\label{propn1}
Define $\Pbnabla_XY$ for $X\in\X(\M), Y\in \Img(P)$ over $\F(\M)$ by
\begin{equation}
\Pbnabla_XY:=\Pnabla_{P(X)}Y+P([X-P(X),Y]).
\end{equation}
Then $\Pbnabla_XY$ has the covariant derivative properties. \qed
\end{propn}

Note that if $\Img(P)$ is a distribution of constant rank, and so is a vector sub-bundle of the tangent bundle, and if indeed $\Pbnabla_XY \in \Img(P)$ for all $X \in \X(M)$ and $Y \in \Img(P)$, then $\Pbnabla$ is just an ordinary linear connection on the vector bundle $\Img(P) \to \M$.


Now suppose that $T\M$ has a direct sum decomposition into $N$ distributions (not necessarily of the same rank) with corresponding unique projectors $P_A, A=1,\dots,N$.
If each distribution has a covariant derivative $\nabla^A$ with extension ${\bar\nabla}^A$ as defined in proposition \ref{propn1} then:
\begin{thm} \label{thm1} For $X,Y\in\X(\M)$
\[\nabla_XY:=\sum_{A=1}^N {\bar\nabla}^A_X(P_A(Y))\]
has the covariant derivative properties over $\F(\M)$. \qed
\end{thm}

That is to say, $\nabla$ is just an affine connection on the manifold $\M$.


The following corollary will be useful in the direct sum situation of theorem \ref{thm1}. The proof can be found in section 2 of \cite{PPSD21}.

\begin{cor}\label{cor2.3}
$\nabla P_B=0$ for all $B\in \{1,\dots,N\}$ if and only if $\nabla^B_XY\in\Img(P_B)$ for all $X,Y\in\Img(P_B)$ and for each $B\in\{1,\dots,N\}$. \qed
\end{cor}




Proposition \ref{propn1} and theorem \ref{thm1} provide the building blocks for the covariant derivatives we will construct in the remaining sections.

\section{The construction of covariant derivatives II}

Suppose that we have a pair of disjoint distributions $K,L$ of equal constant rank on $\M.$  (We don't assume that $T\M=K\oplus L.$) Let $\bar K, \bar L$ denote fixed complements of $K,L$ in $T\M$ with $K \subset \bar{L}$ and $L \subset \bar{K}$, defining projectors $P_K$ and $P_L$. Now assume the existence of a {\em $K$-endomorphism}, a $K$-valued 1-form $S$ on $T\M$ with $\Ker(S)=\bar L, \Img(S)=K$. This determines a unique $L$-valued 1-form ({\em $L$-endomorphism}) $Q$ with $\Ker(Q)=\bar K$, $\Img(Q)=L$ and $Q\circ S=P_L$, $S\circ Q=P_K.$

\begin{thm}\label{thm2}
The derivations $\nabla^K, \nabla^L$ on $K$ and $L$ respectively defined by
\[
\nabla^K_XY:=S([X,Q(Y)]),\ X,Y\in K;\quad \nabla^L_XY:=Q([X,S(Y)]),\ X,Y\in L
\]
\end{thm}
have the covariant derivative properties over $\F(M)$.

\begin{proof}
It suffices to show that $\nabla^K$ has the properties. The $\R$-linearity property of $\nabla^K_XY$ is obvious.
For $X,Y\in K$ and $f\in \F(\M)$
\begin{align*}
\nabla^K_{fX}Y&=S([fX,Q(Y)])=S(f[X,Q(Y)]-Q(Y)(f)X)\\
&=fS([X,Q(Y)])\ \text{since}\ K \subset \bar L = \Ker(S) \\
&=f\nabla^K_XY
\end{align*}
and
\begin{align*}
\nabla^K_X(fY&)=S([X,fQ(Y)])=S(X(f)Q(Y)+f[X,Q(Y)])\\
&=X(f)Y+fS([X,Q(Y)])\ \text{since}\ P_K=S\circ Q\\
&=X(f)Y+f\nabla^K_XY
\end{align*}
\end{proof}

The following generalisation will be useful in the case of Ehresmann connections.

Suppose $\M$ has disjoint distributions $L_1,\dots,L_N$ and $K$ of equal constant rank. Let $\bar K, \bar L_A, A=1,\dots N,$ denote fixed complements of $K,L_A$ in $T\M$ with $L_A \subset \bar{K}$ and $L_B,K \subset \bar{L}_A$ for $B \ne A$, and denote the corresponding projectors by $P_K, P_{L_A}.$ What we have in mind is that $L:=L_1\oplus\dots\oplus L_N$ has some geometric meaning such as an $N$-fold splitting of the vertical or horizontal distribution of an Ehresmann connection.

Further, assume the existence of $N$ $K$-endomorphisms, $S_1,\dots,S_N,$ on $T\M$ with $\Ker(S_A)=\bar L_A$ and $\Img(S_A)=K$. These determine $M$ unique $L$-endomorphisms $Q_A$ with $\Ker(Q_A)=\bar K, \Img(Q_A)=L_A$ and $Q_A\circ S_A=P_{L_A}$ and $S_A\circ Q_A=P_K;$ notice that $S_A\circ Q_B=0, A\neq B.$ Finally, set $S:=
\sum_AS_A / \sqrt{N}$ and $Q:=
\sum_AQ_A / \sqrt{N}$ and note that $S\circ Q=P_K.$

\begin{propn}\label{propn2}
The derivations  $\nabla^K, \nabla^{L_A}$ on $K$ and $L_A$ respectively defined by
\begin{align*}
\nabla^K_XY&:=S([X,Q(Y)]),\ X,Y\in K\\
\nabla^{L_A}_XY&:=Q_A([X,S_A(Y)]),\ X,Y\in L_A,
\end{align*}
have the covariant derivative properties over $\F(M)$. \qed
\end{propn}
Note that $\nabla^K_XY:=S_A([X,Q_A(Y)])$ for any $A$ also has the desired property but is less canonical.

\section{Special Ehresmann Connections}

Let  $\pi:E\to B$ be a fibred manifold with an Ehresmann connection: $V$ is the vertical sub-bundle of $TE$ comprised of tangent spaces to the fibres $E_x$, and $H$ is the horizontal sub-bundle of $TE$ complementary to $V.$ Remembering that pointwise $TE$ is a direct sum of $V$ and $H$,  if both $V$ and $H$ have covariant derivatives $\nabla^V, \nabla^H$ then a covariant derivative on $E$ can be constructed using theorem \ref{thm1}. Of course this begs the question of the construction of these separate derivatives; in geometric cases $V$ will have a covariant derivative by virtue of the geometry of the fibres $E_x$ but for $H$ this requires additional structure.

The significant question is whether such covariant derivatives on the total space $E$ are restricted to the vector bundle case.  We shall give examples in the next section.

If $\dim(V)=\dim(H),$ that is, the dimensions of the tangent spaces to $E_x$ and $B$ coincide (for example when $E=TB$) then the existence of a vertical endomorphism is sufficient to construct a covariant derivative (see theorem \ref{thm3} below). In the case where $\dim(V)<\dim(H)$ a vertical endomorphism can still be used if $H=H_1\oplus H_2$ with $\dim(H_1)=\dim(V)$ and where $H_2$ is equipped with a covariant derivative. While this sounds a little artificial, it is exactly the case which arises in the geometry of time-dependent second order ordinary differential equations. Indeed, in that case the vertical endomorphism is central and $V, H_1$ and $H_2$ naturally appear as eigenspaces of a related endomorphism, see \cite{PPSD21, KP08}. We shall treat the simpler, time-independent case in example \ref{SODE} below.

The next result follows directly from proposition \ref{propn1} and theorems \ref{thm1},\ \ref{thm2} with $K,L=V,H$.

\begin{thm}\label{thm3}
Let $\pi:E\to B$ be a fibred manifold with fibres $E_x$ and base $B$ having the same dimension. Assume an Ehresmann connection on $E$ with vertical and horizontal distributions $V, H$ and corresponding projectors $P_V, P_H$. Further suppose that $S$ and $Q$ are vertical and horizontal endomorphisms (so that $\Img(S)=V=\Ker(S)$ and $\Img(Q)=H=\Ker(Q)).$ Then $E$ has covariant derivative
\begin{align}
\nabla_XY&:=\bnabla^V_X(P_V(Y))+\bnabla^H_X(P_H(Y))\notag\\
&:=S([P_V(X),Q(Y)])+P_V([P_H(X),P_V(Y)])\\
&\quad +Q([P_H(X),S(Y)])+P_H([P_V(X),P_H(Y)]).\notag
\end{align}
\qed
\end{thm}

Recalling that the curvature of an Ehresmann connection is $R(X,Y)=P_V([P_H(X),P_H(Y)])$ it is straightforward to see that the torsion of this covariant derivative also carries this information: $P_V\circ T(P_H(X),P_H(Y))=R(X,Y)$. This will appear again in the examples.


Another straightforward case is when $\dim(H)$ is an integer multiple of $\dim(V)$ or vice-versa.

\begin{thm}\label{thm4}
Let $\pi:E\to B$ be a fibred manifold with fibres $E_x$ and base $B$ with $\dim(B)$ being an integer multiple $N$ of $\dim(E_x)$.
Assume an Ehresmann connection on $E$ with vertical subbundle $V$ and horizontal subbundle $H=H_1\oplus\dots\oplus H_N$ with $\dim(H_A)=\dim(V)$ for all $A.$
Let the corresponding projectors, vertical and horizontal endomorphisms be $P_V, P_{H_A}, S_A, Q_A$ with $S:=
\sum_AS_A / \sqrt{N}$, $Q:=
\sum_AQ_A / \sqrt{N}$.
Then $E$ has a covariant derivative
\begin{align}
\nabla_XY&:=\bnabla^V_X(P_V(Y))+\sum_A\bnabla^A_X(P_A(Y))\notag\\
&:=S([P_V(X),Q(Y)])+\sum_AP_V([P_A(X),P_V(Y)])\\
&\quad +\sum_A\left( Q_A([P_A(X),S_A(Y)])+P_A([X-P_A(X),P_A(Y)]) \right).\notag
\end{align}
A similar result holds, with relabelling, when $\dim(E_x)$ is an integer multiple of $\dim(B)$.
\end{thm}
The proof when $\dim(H) = N \dim(E_x)$ follows from proposition \ref{propn1}, theorem \ref{thm1} and proposition \ref{propn2} with $K=V$ and $L_A=H_A$.

The case where $V$ is an $N$-fold direct sum with $\dim(V_A)=\dim(H)$ requires $L_A=V_A$ and $K=H$, and also the switching of labels $Q$, $S$ so that $S_A$ labels the vertical endomorphisms and $Q_A$ labels the horizontal ones as usual.
\qed

As a result of corollary \ref{cor2.3} the covariant derivatives of theorems \ref{thm3} and \ref{thm4} have the properties $\nabla P=0$ for all the projectors.
In addition, in the case of the covariant derivative of theorem \ref{thm3} $\nabla S=0=\nabla Q.$ This is not so in the case of theorem \ref{thm4}.

\section{Applications and Examples}

In this section we give a variety of examples of our construction. We start with some general ones and then consider, in particular, examples on a tangent bundle, and finally on a frame bundle.

\subsection{General examples}

\begin{xmpl} Trivial bundle\newline
Let $\pi : \R^3 \to \R^2$ be the trivial bundle with base coordinates $x$, $y$ and fibre coordinate $\theta$, and take as a horizontal distribution
\[
H(E)=\Sp\left\{H_1:=\pd{}{x}+\cos\theta\pd{}{\theta},\ H_2:=\pd{}{y}+\sin\theta\pd{}{\theta}\right\}
\]
so that $H$ satisfies the conditions for an Ehresmann connection.

The vertical and horizontal endomorphisms for use in theorem \ref{thm4} will be $S_1:=\d x\otimes V, S_2:=\d y\otimes V$ and $Q_1:=\psi\otimes H_1, Q_2:=\psi\otimes H_2$ with $V := \partial / \partial \theta$ and $\psi:=\d\theta - \cos\theta \, \d x + \sin\theta \, \d y$. The nonzero components of $\nabla$ will then be
\[
\nabla_{H_1}H_1 = \sin\theta \, H_1, \quad
\nabla_{H_2}H_2 = -\cos\theta \, H_2, \quad
\nabla_{H_1}V = \sin\theta \, V, \quad
\nabla_{H_2}V = -\cos\theta \, V.
\]
\end{xmpl}

\begin{xmpl} The Hopf bundle\newline
A more complex application of theorem \ref{thm4} is to the Hopf bundle $\pi : S^3 \to S^2$. If we regard $S^3$ as the submanifold of $\R^4$ given by $x^2 + y^2 + z^2 + w^2 = 1$ then one choice for the definition of $\pi$ is given by
\[
\pi(x,y,z,w) = \bigl( x^2 + y^2 - z^2 - w^2, 2(xw + yz), 2(yw - xz) \bigr) \in S^2 \subset \R^3 \, .
\]
Now consider the infinitesimal rotations in the six coordinate $2$-planes of $\R^4$,
\begin{gather*}
\theta = y \vf{x} - x \vf{y} \, , \qquad \phi = z \vf{x} - x \vf{z} \, , \qquad \psi = w \vf{x} - x \vf{w} \, , \\
\xi = w \vf{z} - z \vf{w} \, , \qquad \eta = y \vf{w} - w \vf{y} \, , \qquad \zeta = z \vf{y} - y \vf{z}
\end{gather*}
so that these restrict to vector fields tangent to $S^3$, and put
\[
V = \theta - \xi \, , \qquad \Lambda = \phi - \eta \, , \qquad \Sigma = \psi - \zeta
\]
as vector fields on $S^3$. We may check that $(\Lambda, \Sigma, V)$ form a global frame for $S^3$ defining an absolute parallelism, and that $V$ is vertical with respect to $\pi$, so that $(\Lambda, \Sigma)$ span the horizontal sub-bundle of an Ehresmann connection on the Hopf bundle with
\[
[\Sigma, \Lambda] = 2V \, , \qquad [\Lambda, V] = 2\Sigma \, , \qquad [V, \Sigma] = 2\Lambda \, .
\]
Writing $H_1 = \langle \Lambda \rangle$ and $H_2 = \langle \Sigma \rangle$, we can apply theorem \ref{thm4} to obtain a covariant derivative $\nabla$ on $S^3$. It is straightforward to check that if $X,Y \in \{ \Lambda, \Sigma, V \}$ then $\nabla_X Y = 0$. This linear connection is not symmetric, but if we subtract half the torsion to obtain a symmetric connection with the same geodesics then the result is just the Levi--Civita connection of the standard metric on $S^3$ induced by its embedding in $\R^4$.
\end{xmpl}

\subsection{Tangent bundle connections}
The next examples relate to various tangent bundle connections; we will uncover the central role of torsion in relating affine, \SODE\ and spray connections to the general, nonlinear case.  
We begin with the tangent bundle connection induced by an affine connection on the base manifold.

\begin{xmpl} Affine connection with torsion\label{Affcontor}\newline
Consider an affine connection with torsion on a manifold $M^n$, with connection coefficients $\Gamma^a_{bc}$ in local coordinates $x^a.$ This induces an $n$-dimensional  horizontal distribution on the tangent bundle $TM$ (with local adapted coordinates $(x^a, u^a)$) spanned by
\[
H_a:=\pd{}{x^a}-\Gamma^c_{ab}u^b\pd{}{u^c}.
\]
With vertical distribution $\Sp\{V_c:=\partial / \partial u^c\},$ the basis dual to $\{H_a,V_b\}$ is
\[
\{\d x^a, \phi^b:=\d u^b+ \Gamma^b_{cd}u^d\d x^c\}.
\]
The non-zero bracket relations are
\[[H_a,V_b]=\Gamma^c_{ab}V_c,\ [H_a,H_b]=\left[\left(\pd{\Gamma^c_{ad}}{x^b}-\pd{\Gamma^c_{bd}}{x^a}\right)+(\Gamma^e_{ad}\Gamma^c_{be}-\Gamma^e_{bd}\Gamma^c_{ae})\right]u^dV_c
\]
Notice that $[H_a,H_b]=-R^c_{dab}u^dV_c$ where $R^c_{dab}$ are the components of the Riemann curvature of the connection on $M^n$.

It is straightforward to verify that the (1,1) tensors
\[
S=\d x^a\otimes V_a\ \text{and}\ Q=\phi^a\otimes H_a
\]
are intrinsic (tensorial). Using these in the definitions of $\nabla^V,\nabla^H$ from theorem \ref{thm2}, we find
\[
\nabla^V_{V_a}V_b=0,\ \nabla^H_{H_a}H_b=\Gamma^c_{ab}H_c.
\]

Constructing the covariant derivative $\nabla$ of theorems \ref{thm1} and \ref{thm3} with $K,L=V,H$ gives
\[
\nabla_{V_a}V_b=0,\ \nabla_{V_a}H_b=0,\ \nabla_{H_a}V_b=\Gamma^c_{ab}V_c,\ \nabla_{H_a}H_b=\Gamma^c_{ab}H_c.
\]
A straightforward calculation gives the only non-zero torsion components as
\[
T(H_a,H_b)=(\Gamma^c_{ab}-\Gamma^c_{ba})H_c -\left[\left(\pd{\Gamma^c_{ad}}{x^b}-\pd{\Gamma^c_{bd}}{x^a}\right)+(\Gamma^e_{ad}\Gamma^c_{be}-\Gamma^e_{bd}\Gamma^c_{ae})\right]u^dV_c
\]
which capture, as horizontal and vertical lifts respectively, the torsion and curvature of the original affine connection on $M.$ Note also that the Ehresmann curvature is the vertical part of the torsion (up to a sign).
\end{xmpl}

\begin{xmpl} Nonlinear connection on $TM$ \label{nonlinTM}\newline
In the previous example the horizontal distribution on $M$ arising from an affine connection was just an Ehresmann connection on the tangent bundle, linear with respect to the vector bundle structure. More generally, we can consider an Ehresmann connection on $TM$ which need not be linear, and apply \ref{thm3} directly with $K,L=V,H$. As in that earlier example, we shall take coordinates $(x^a,u^b)$ on $TM$. Put $V=\Sp\{V_a:=\partial / \partial u^a \}$ and
\[
H=\Sp\left\{H_a:=\pd{}{x^a}-\Gamma^b_a\pd{}{u^b}\right\}
\]
so that
$S=\d x^a\otimes V_a$ and $Q=\phi^a\otimes H_a,$ where $\{\d x^a, \phi^b:=\d u^b+ \Gamma^b_c\d x^c\}$ is dual to $\{H_a,V_b\}$. The results appear similar to those in the affine connection case with $\Gamma^a_b$ replacing $\Gamma^a_{bc}u^c$. We get

\[
\nabla^V_{V_a}V_b=S([V_a,H_b])=0,\ \nabla^H_{H_a}H_b=Q([H_a,V_b])=V_b(\Gamma^c_a)H_c
\]

and
\[
\nabla_{V_a}V_b=0,\ \nabla_{V_a}H_b=0,\ \nabla_{H_a}V_b=V_b(\Gamma^c_a)V_c,\ \nabla_{H_a}H_b=V_b(\Gamma^c_a)H_c.
\]
More interesting is that the torsion acts non-trivially only on $H$ with
\[
T(H_a,H_b)=Q([H_a,V_b]-[H_b,V_a])-[H_a,H_b];
\]
this means that the horizontal component is
\[
P_H(T(H_a,H_b))=\left(V_b(\Gamma^c_a)-V_a(\Gamma^c_b)\right)H_c
\]
and the vertical component is the Ehresmann curvature $R(H_a,H_b).$ This horizontal component is zero if and only if there are local functions
$f^c$ such that $\Gamma^c_a = V_a \bigl( - \onehalf f^c \bigr)$. The reason for the numerical coefficient, and a global explanation of the phenomenon, may be seen in the next example.
\end{xmpl}

\begin{xmpl} \SODE\ connection on $TM$ \label{SODE}\newline
Again we take $E:=TM$ with standard vertical endomorphism $S=\d x^a\otimes V_a.$ A second order differential equation field (\SODE) is a vector field $\Gamma$ on $TM$ satisfying $S(\Gamma) = \Delta$, where $\Delta$ is the canonical dilation field on the vector bundle $\tau_M : TM \to M$ \cite{C83}. The reason for the name is that such a vector field defines a family of curves on $M$ satisfying second order equations of the form $\ddot{x}^a = f^a(x,\dot{x})$.

Every \SODE\ gives rise to a nonlinear connection on $TM$ with horizontal projector $\PH = \onehalf \bigl( I - \lie{\Gamma} S \bigr)$; this is called the \SODE\ connection of $\Gamma$. In coordinates, if
\[
\Gamma = u^a \fpd{}{x^a} + f^a \fpd{}{u^a}
\]
then
\[
\PH = dx^a \otimes \biggl( \fpd{}{x^a} - \Upsilon^b_a \fpd{}{u^a} \biggr) \, , \qquad
\Upsilon^b_a = - \onehalf \fpd{f^b}{u^a} \, .
\]
We might therefore ask when a general nonlinear connection on $TM$ was derived in this way from a \SODE, and then we can use theorem~\ref{thm3} as in the previous example to construct an affine version. So suppose such a nonlinear connection has horizontal vector fields
\[
H_a = \fpd{}{x^a} - \Gamma^b_a \fpd{}{u^b} \, .
\]
Regard the connection as a splitting of the short exact sequence of vector bundles over $TM$
\[
\begin{tikzcd}
0 \arrow{r} & V\tau_M \arrow{r} & TTM \arrow{r} & \tau_M^* TM \arrow{r} \arrow[bend right, swap]{l}{h} & 0
\end{tikzcd}
\]
where, if $v, w \in T_p M$, the splitting map $h$ gives the horizontal lift of $v$ to $h(v,w) \in T_w TM$, so that
\[
v = v^a \vfe{x^a}{p} \, , \qquad h(v,w) = v^a H_a|_w = v^a \biggl( \vfe{x^a}{w} - \Gamma^b_a(w) \vfe{u^b}{w} \biggr) \, .
\]
Define a vector field $\Gamma$ on $TM$ by $\Gamma_v = h(v,v) \in T_v TM$, so that
\[
\Gamma_v = v^a \biggl( \vfe{x^a}{v} - \Gamma^b_a(v) \vfe{u^b}{v} \biggr)
= u^a(v) \biggl( \vfe{x^a}{v} - \Gamma^b_a(v) \vfe{u^b}{v} \biggr)
\]
and therefore that
\[
\Gamma = u^a H_a = u^a \vf{x^a} - u^a \Gamma^b_a \vf{u^b} \, .
\]
We see that $\Gamma$ is a \SODE\ with $f^b := - u^a \Gamma^b_a$, giving rise to a \SODE\ connection. Put
\[
\Upsilon^b_c = - \tfrac{1}{2} \pd{f^b}{u^c} = \tfrac{1}{2} \vf{u^c} (u^a \Gamma^b_a)
= \tfrac{1}{2} \bigl( \Gamma^b_c + u^a V_c(\Gamma^b_a) \bigr) \, .
\]
Now suppose that the original nonlinear connection is `symmetric' (that is, that the torsion of $\nabla$ is vertical, as in the previous example) and that the connection coefficients $\Gamma^b_a$ are homogeneous (in other words that $\Delta(\Gamma^b_a) = \Gamma^b_a$, where $\Delta = u^c V_c$). We see that in these circumstances
\[
u^a V_c(\Gamma^b_a) = u^a V_a(\Gamma^b_c) = \Delta(\Gamma^b_c) = \Gamma^b_c
\]
so that $\Gamma^b_c = \Upsilon^b_c$ and the two connections are identical. Note, incidentally, that
\[
\Delta(f^b) = - u^a \Gamma^b_a - u^a \Delta(\Gamma^b_a)
\]
so that $\Gamma$ is a spray, with $\Delta(f^b) = 2f^b$, precisely when the connection coefficients are homogeneous.

Now return to the affine connection on $TM$. We find that
\[
\nabla_{H_a} \Gamma = \bigl( \Delta(\Gamma^b_a) - \Gamma^b_a \bigr) H_b \, , \qquad
\nabla_{V_a} \Gamma = H_a
\]
so that the connection coefficients are homogeneous precisely when $\nabla_{X^H} \Gamma = 0$ for any vector field $X$ on $M$, and similarly
\[
\nabla_{H_a} \Delta = \bigl( \Delta(\Gamma^b_a) - \Gamma^b_a \bigr) V_b \, , \qquad
\nabla_{V_a} \Delta = \nabla_{V_a} (u^b V_b) = V_a
\]
so that the connection coefficients are also homogeneous precisely when $\nabla_{X^H} \Delta = 0$ for any vector field $X$ on $M$.

\begin{propn}
If $\nabla$ is the affine connection on $TM$ obtained from a nonlinear connection as in example~\ref{nonlinTM} then
\[
\nabla_{X^H} \Delta = 0 \, , \qquad P_H \circ T = 0
\]
is a sufficient condition for the nonlinear connection to be a \SODE\ connection, and it will then be the \SODE\ connection of a spray. \qed
\end{propn}
\end{xmpl}

\subsection{Frame Bundle Connections}

We begin with the frame bundle $\pi : F \to M$ with $\dim M = n$, where each fibre $F_x$ is the set of frames (bases) for the tangent space to $M$ at $x$ (see \cite{CP86, CCL99}). The general linear group $GL(n,\R)$ acts freely and transitively on each fibre, by matrix multiplication, conventionally on the right. The fibre can be identified locally with $GL(n,\R)$ using a fixed co-ordinate basis, but generally there is no way of globally fixing the identity of $GL(n,\R)$ on each fibre. The dimension of $F$, and therefore the dimension of any tangent space $T_p F$, is $n^2 + n$. The {\em vertical subspace} of $T_pF$ has dimension $n^2$ and is the tangent space at $p$ to the fibre $F_{\pi(p)}$; it may be identified with $M_{n\times n}(\R)$. Any horizontal subspace will have dimension $n$ and we will only be able to apply theorem \ref{thm3} if we can find an $n$-fold decomposition of the vertical subspace. To that end we present, without proof, this special case of the Primary Decomposition Theorem of linear algebra.

\begin{propn}\label{perm'n prop}
Let $A$ be an $n \times n$ matrix obtained from the identity matrix $I_n$  by permuting its rows using a permutation with a single orbit (in other words, an $n$-cycle),
and define $T_A$ to be the linear operator on $M_{n\times n}(\R)$ with $T_A(B):=AB.$ Then
\[M_{n\times n}(\R)=W^1\oplus\dots\oplus W^n\]
with $W^1,\dots,W^n$ being invariant subspaces of $T_A,$ each of dimension $n$ and corresponding to an $n^{th}$ root of unity linear factor in the minimal polynomial over $\mathbb{C}$ of $T_A$.
Moreover, all such linear operators have the same invariant subspaces, where $W^k$ contains matrices whose nonzero entries are all in column $k$.
\end{propn}

Now, with any such permutation of the identity in mind, apply the result of proposition \ref{perm'n prop} to $M_{n\times n}(\R)$ identified as the Lie algebra of $GL(n,\R).$ Using a basis associated with the invariant subspaces of the decomposition, create the usual fundamental vector fields for each of the corresponding $n$-parameter subgroups of $GL(n,\R)$. Because these vector fields are tangent to the fibres they are vertical, moreover they span the tangent spaces at each point on a fibre because the right action of the subgroups act transitively on the fibres. In this way we create an $n$-fold splitting of the vertical distribution of $\pi : F \to M$.

The final ingredient required is the horizontal distribution. We shall follow \cite{CCL99} and choose an affine connection on $M,$ but we could use any of the connections discussed in the previous section.

We begin with bases $\{W_1^A,\dots,W_n^A\}$
for the invariant subspaces $W^A$ of $M_{n\times n}(\R),$ along with coordinates $(x^i)$  and co-ordinate frame $\{\partial / \partial x^1,\dots,\partial / \partial x^n\}$ for open $U\subset M.$ The components of the affine connection relative to this frame are $\Gamma^i_{jk}$ as in example~\ref{Affcontor}. Now a smooth pointwise change of frame  to any arbitrary basis $\{b_1,\dots,b_n\}$ is represented by a non-singular matrix $(b)^i_j$ where $b^i_j\in\mathfrak{F}(U)$ and
\[b_j=b^i_j\vf{x^i}.\]

In this way we obtain local co-ordinates $(x^i,b^i_j)$ for the frame bundle $F$. Now we wish to utilise the fixed invariant subspaces $W^A,$ to this end we write $(b)^i_j=w^a_B(W_a^B)^i_j$ where $w^a_B\in\mathfrak{F}(U),$ producing coordinates $(x^i,w^a_B)$ for $\pi^{-1}(U)\subset F$. The $n^2$-dimensional vertical subspaces are locally generated by $V^A_b:=\partial / \partial w^b_A$ where the upper index on $V^A_b$ labels the invariant subspace $W^A.$

Finally, we have local generators for the $n$-dimensional horizontal subspaces of the tangent spaces to $F$:
\[H_i:=\vf{x^i}-\Gamma^k_{ij}w^j_A\vf{w^k_A}.\]

The non-zero bracket relations are
\[[H_a,V^A_b]=\Gamma^c_{ab}V^A_c,\ [H_a,H_b]=\left[\left(\pd{\Gamma^c_{ad}}{x^b}-\pd{\Gamma^c_{bd}}{x^a}\right)+(\Gamma^e_{ad}\Gamma^c_{be}
-\Gamma^e_{bd}\Gamma^c_{ae})\right]w_A^dV^A_c.
\]

The induced vertical and horizontal endomorphisms of the `flipped' version of theorem \ref{thm4} are
\[S^A=\d x^b\otimes V^A_b\ \text{and}\ Q_A=\phi^b_A\otimes H_b\]
where $\phi^b_A:=\d w^b_A+\Gamma^b_{cd}w^d_A\d x^c.$ Using these in the definitions of $\nabla^V,\nabla^H$ from theorem \ref{thm2}, we find
\[
\nabla^{V^A}_{V^A_a}V^A_b=0,\ \nabla^H_{H_a}H_b=\Gamma^c_{ab}H_c.
\]
Calculating the resulting covariant derivative produces 
\[
\nabla_{V^A_a}V^B_b=0,\ \nabla_{V^A_a}H_b=0,\ \nabla_{H_a}V^B_b=\Gamma^c_{ab}V^B_c,\ \nabla_{H_a}H_b=\Gamma^c_{ab}H_c,
\]
which is close to the tangent bundle case, as is the only non-zero component of the torsion

\[
T(H_a,H_b)=(\Gamma^c_{ab}-\Gamma^c_{ba})H_c -\left[\left(\pd{\Gamma^c_{ad}}{x^b}-\pd{\Gamma^c_{bd}}{x^a}\right)+(\Gamma^e_{ad}\Gamma^c_{be}-\Gamma^e_{bd}\Gamma^c_{ae})\right]w^d_AV^A_c,
\]
vindicating this formulation of the frame bundle connection.

\section*{Acknowledgement}
The authors thank a referee for useful suggestions regarding the layout and content of this paper.

G.E.\,Prince \\
Department of Mathematics and Statistics, La Trobe University,\\
Victoria 3086, Australia \\
Email: \url{g.prince@latrobe.edu.au}

D.J.\,Saunders \\
Lepage Research Institute,\\
17. novembra 1, 081 16 Pre\v{s}ov, Slovakia \\
Email: \url{david@symplectic.email}

\end{document}